\documentclass[mathpazo]{ata}
\newtheorem{thm}{Theorem}[section]

\newtheorem{defi}{Definition}[section]
\newtheorem{deduction}{Corollary}[section]
\newtheorem{rem}{Remark}[section]
\usepackage{comment}
\usepackage{CJKutf8}
\usepackage{graphicx}
\usepackage{subfigure}
\usepackage{float}
\usepackage{booktabs} 
\usepackage{mathrsfs}
\usepackage{amssymb}
\graphicspath{{figures/},{pics/}}
\usepackage{xcolor}
\begin{document}
\begin{CJK}{UTF8}{gbsn}
\title{Estimation of sparse polynomial approximation error to continuous function}


 \author[Feng and Zhang]{Renzhong Feng\affil{1}, Bowen Zhang\affil{1}\comma\corrauth}
 \address{\affilnum{1} School of Mathematical Sciences \& LMIB of the Ministry of Education, Beihang University, 37 Xueyuan Road, Haidian District, Beijing, 100191, China}
 \emails{{\tt fengrz@buaa.edu.cn} (R.~Feng), {\tt zhangbowensy@163.com} (B.~Zhang)}

\begin{abstract}
	
The sparse polynomial approximation of continuous functions has emerged as a prominent area of interest in function approximation theory in recent years. A key challenge within this domain is the accurate estimation of approximation errors. This paper focuses on continuous functions, characterizing their sampled values as a combination of the values of their best approximation polynomials within a finite-dimensional polynomial space and the associated remainder terms. Consequently, the sampled values of a function can be interpreted as noisy samples of the values of its best approximation polynomial, with the noise equivalent to the remainder term's values at those points. By selecting a uniformly bounded orthonormal polynomial system as the basis for this finite-dimensional space, it becomes feasible to formulate noise constraint inequalities of the form $||A\mathbf{c}-\mathbf{y}||_2\leqslant \Theta\sqrt{m}E_N$ and $l_1$-minimization problems or their weighted $l_1$-minimization variants. This paper provides estimations for the approximation error of the sparse polynomial derived from the $l_1$-minimization method, both in $L_\infty$-norm and $L_2$-norm, characterizing the error in terms of the quasi-norm of the sampled function or its best uniform approximation polynomial, the sparsity, and the best approximation error. When the smoothness of the sampled function is known, the approximation error can be further simplified to depend on the quasi-norm of the sampled function or its best uniform approximation polynomial and the sparsity. The analysis reveals that if the sampled function is a sparse polynomial from a finite-dimensional space, it can be reconstructed exactly. Moreover, it is observed that the smoother the sampled function, the fewer degrees of the sparse polynomial are required to attain a given approximation accuracy. The paper also extends this analysis to estimate the $L_2$-norm approximation error for the sparse polynomial obtained \textit{via} the weighted $l_1$-minimization method, noting that in this context, the orthonormal polynomial system does not need to be uniformly bounded for the conclusions to hold.
\end{abstract}

 \clc{O241.5} 

\maketitle

\section{Introduction}
Ben Adcock highlighted in \cite{Ben2017} that numerous smooth functions exhibit sparse coefficients when expanded in terms of orthogonal polynomial bases. This suggests that smooth functions can be effectively approximated by sparse polynomials within orthogonal bases. Concurrently, advancements in compressive sensing (CS) have demonstrated the feasibility of recovering high-dimensional vectors with specific low-dimensional structures from a substantially reduced set of linear measurements \cite{EJ2006,DL2006}. Consequently, there has been a growing interest in employing compressive sensing techniques to compute sparse polynomials \cite{AH2011,JA2015,LK2012,AT2015,JJA2014,HR2011,HR2012,HR2016,LLD2012,XG2013,ARS2022,RAMZ2023}. The standard procedure involves initially conducting a generalized Fourier series expansion of the function to be approximated within an orthogonal polynomial system, then truncating the expansion after the first $N$ terms. The function values at the sampling points are then treated as noisy samples of the truncated series, with the noise equivalent to the values of the remainder function at those points. This approach naturally gives rise to an $l_1$-minimization problem with an inequality constraint, where the solution vector is used to approximate the coefficients of the first $N$ terms, yielding a sparse polynomial approximation of the sampled function. As mentioned in \cite{Ben2017}, the approximation error of the sparse polynomial hinges on the selection of $\delta$ in the inequality constraint, representing the magnitude of the remainder values. Thus, accurately estimating the norm of the tail of the expansion is vital for ensuring the precision of the outcomes. 

A challenge arises here: this tail error is generally unknown. Cross-validation techniques \cite{AH2011,JJA2014,XG2013} are often employed to provide practical estimates for $\delta$, but they are computationally intensive and inefficient in terms of data utilization. This paper employs the $l_1$-minimization problem to provide sparse polynomial approximations for general continuous functions, representing the remainder as the discrepancy between the function to be approximated and its best approximation in the polynomial space of degree $N$. This allows for the convenient use of approximation algorithms, such as the Remez algorithm, to more accurately estimate the norm of the remainder. The approximation error is then expressed as the sum of the best $s$-term approximation error of the first $N$ terms and the best approximation error.
 Using Stechkin's inequality and Jackson's theorem, the approximation error can be further expressed solely in terms of the quasi-norm of the function to be approximated or its best uniform approximation polynomial and the sparsity $s$. These error expressions reveal that the smoother the sampled function, the fewer degrees of the sparse polynomial are required to attain the same level of approximation accuracy.

Regarding the weighted $l_1$-minimization problem, it is distinct from the standard $l_1$-minimization problem in its utilization of a weighted $l_1$-norm. Its advantage lies in effectively mitigating the curse of dimensionality that often occurs with $l_1$-minimization in multivariate situations, and it also has a lower requirement for the orthogonal polynomial system in sparse polynomial recovery. Within this domain, an analogous recovery theory to the unweighted case has been established, successfully extending numerous concepts and associated properties from the unweighted context to the weighted scenario, as detailed in \cite{HR2016}. Building on this, the weighted $l_1$-minimization method has yielded various results from different perspectives in the field of sparse recovery: Adcock and colleagues introduced the concept of hyperbolic cross index sets, deriving a lower bound on the required sampling number for the weighted $l_1$-minimization method where the dimension $d$ no longer appears in the exponent \cite{Ben2018}, thus alleviating the curse of dimensionality we mentioned. Moreover, Adcock and his colleagues have devised a method for selecting weights that ensures the approximation results in noiseless sampling exhibit interpolation properties. This approach ensures that the truncation strategy is independent of the function being approximated and can be applied to any selected sampling points, as discussed in \cite{Ben2017}. These research outcomes undoubtedly highlight the significant benefits of incorporating weighted considerations.

The remainder of this paper is organized as follows: Section 2 introduces some preliminary knowledge used in this paper, Section 3 offers an estimation of the $L_\infty$-norm of the approximation error in the context of sparse polynomial approximation. Section 4 extends this analysis to provide an estimation of the $L_2$-norm of the approximation error for sparse polynomial approximation and broadens the findings from Section 3 to include the weighted $l_1$-minimization scenario, yielding $L_2$-norm error estimates. Section 5 substantiates the efficacy of these novel methods through numerical experiments, while Section 6 summarizes the conclusions drawn from the study.
 
\section{Preliminary}
\subsection{The Best Uniform Approximation of Continuous Functions and Jackson's Theorem}

Let $\Pi_n(I)$ denotes the set of all polynomials defined on $I=[a,b]$ with degree not exceeding $n$. For $f(x)\in C[a,b]$, let $||f||_\infty:=\mathop{\max}\limits_{x\in [a,b]}|f(x)|$. The quantity

$$\Delta(p):=\mathop{\max}\limits_{x\in [a,b]}|f(x)-p(x)|$$	
refers to the deviation of the polynomial $p(x)$ from the function $f(x)$. The quantity
$$E_n=E_n(f):=\mathop{\inf}\limits_{p\in\Pi_n(I)}\{\Delta(p)\}$$
refers to as the minimum deviation of $\Pi_n(I)$ from $f(x)$, or the best uniform approximation of $f(x)$ to $\Pi_n(I)$.

\begin{rem}
	It is evident that $E_n\geqslant 0$ and $E_0\geqslant E_1\geqslant E_2\geqslant\cdots$, and moreover, $\mathop{\lim}\limits_{n\rightarrow\infty}E_n=0.$
\end{rem}

\begin{defi}
	Let $f(x)$ be a continuous function defined on interval $[a,b]$, and let $\delta$ be any positive number. Consider all pairs of points $x,y$ in $[a,b]$ such that $|x-y|\leqslant \delta$. The quantity
    $$\omega(\delta):=\mathop{\sup}\limits_{|x-y|\leqslant \delta}\{|f(x)-f(y)|\}$$
is called the modulus of continuity of function $f(x)$.
\end{defi}

The modulus of continuity is a crucial attribute that reflects the intrinsic structure of a function, as demonstrated by Jackson's theorem \cite{Nata1958} and its corollaries.

\begin{thm}[Jackson's Theorem \cite{Nata1958}]\label{thm1}
	If $E_n$ represents the best uniform approximation of function $f(x)\in C[a,b]$ by polynomials in space $\Pi_n(I)$, then
	$$E_n\leqslant 12\omega\left(\frac{b-a}{2n}\right),$$
	where $\omega(\delta)$ denotes the modulus of continuity of $f(x)$.
\end{thm}

From Theorem \ref{thm1}, we can derive the following two corollaries.

\begin{deduction}[\cite{Nata1958}]\label{coro1}
	If $f(x)\in Lip_M^\alpha(0<\alpha\leqslant 1)$, then	
    $$E_n\leqslant \frac{CM}{n^\alpha},$$
	where $C=12\left(\frac{b-a}{2}\right)^\alpha.$	
\end{deduction}

\begin{deduction}[\cite{Nata1958}]\label{coro2}
	If function $f(x)$ has a bounded derivative $f'(x)$ such that $|f'(x)|\leqslant M_{1}$, then
	$$E_n\leqslant \frac{6(b-a)M_{1}}{n}.$$	
\end{deduction}

\begin{thm}[\cite{Nata1958}]\label{thm2}
	If function $f(x)$ has a bounded $(p+1)$-th derivative $f^{(p+1)}(x)$ such that $|f^{(p+1)}(x)|\leqslant M_{p+1}$, then
	$$E_n\leqslant \frac{C_p(b-a)^{p+1}M_{p+1}}{n^{p+1}},$$
	where the constant $C_p$ depends only on $p$.	
\end{thm}

\subsection{Sparse Polynomial Recovery}

Assume that $I$ is endowed with a probability measure $\nu$. Furthermore, assume that $\psi_1,\dots,\psi_N$ form an orthonormal polynomial system defined on $I$ with respect to $\nu$, that is, for $j,k\in [N]:=\{1,2,\dots,N\}$，
$$\int_\mathbb{D}\psi_{j}(x)\psi_k(x)\nu(x)dx=\left\{
\begin{aligned}
	&0\quad if\hskip 0.5em j\neq k,\\
	&1\quad if\hskip 0.5em j=k.
\end{aligned}
\right.$$

If
$$||\psi_j||_\infty=\mathop{\sup}\limits_{x\in I}|\psi_j(x)|\leqslant K\quad for\hskip 0.5em all \hskip 0.5em j\in [N],$$
where $K\geqslant 1$ is a constant, we refer to the function system $\psi_1,\dots,\psi_N$ as being uniformly bounded. We focus on functions of the following form
$$g(x)=\mathop{\sum}\limits_{k=1}^Nc_k\psi_k(x),\quad x\in\mathbb{D},$$
where the coefficients $c_1,c_2,\dots,c_N\in\mathbb{R}$. If
$$||\mathbf{c}||_0:=|supp(\mathbf{c})|\leqslant s<N,$$
then $\mathbf{c}$ is said to be $s$-sparse. Consequently, function $g(x)$ is called an $s$-sparse polynomial.

Let $x_1,x_2,\dots,x_m$ be a set of sampling points on interval $I$. The corresponding sampled values of function $g(x)$ at these points are given by
$$y_l=g(x_l)=\mathop{\sum}\limits_{k=1}^Nc_k\psi_k(x_l),\quad l=1,2,\dots,m.$$
We define the sampling matrix $A\in\mathbb{R}^{m\times N}$, where each element is
$$A_{l,k}=\psi_k(x_l),\quad l=1,\dots,m,\quad k=1,\dots,N.$$
Consequently, the vector of sampled values $\mathbf{y}=(y_1,\dots,y_m)^T$ can be represented as
$$\mathbf{y}=A\mathbf{c},$$
where $\mathbf{c}=(c_1,\dots,c_N)^T$ denotes the coefficient vector of function $g(x)$.

The $l_p$-norm on $\mathbb{R}^N$ is defined as $||\mathbf{z}||_p:=\left(\mathop{\sum}_{j=1}^N|z_j|^p\right)^{\frac{1}{p}},\quad 1\leqslant p<\infty$，and $||\mathbf{z}||_\infty:=\mathop{\max}\limits_{1\leqslant j\leqslant N}|z_j|$. The best $s$-term approximation error of a vector $\mathbf{z}\in\mathbb{R}^N$ in the $l_p$-space is defined as
$$\sigma_s(\mathbf{z})_p:=\mathop{\inf}\limits_{\mathbf{x}\in\mathbb{R}^N,||\mathbf{x}||_0\leqslant s}||\mathbf{z}-\mathbf{x}||_p,\quad p>0.$$
Stechkin provided an estimate for $\sigma_s(\mathbf{z})_p$ in \cite{MH2011}:
$$\sigma_s(\mathbf{z})_p\leqslant s^{\frac{1}{p}-\frac{1}{q}}||\mathbf{z}||_q,\quad 0<q<p.$$

Furthermore, we introduce the concept of weighted $l_p$-space and weighted $l_p$-norm for $0<p\leqslant 2$. The weighted $l_p$-norm is an extension of the standard $l_p$-norm. When weights are assigned, we consider an index set $\Lambda$ that may be finite or countably infinite, and associate with it a weight vector $\boldsymbol{\omega}=(\omega_j)_{j\in \Lambda}$.  It is customary to assume that each component of the weight vector satisfies $\omega_j\geqslant 1$. The definitions for the weighted $l_p$-norm and the weighted best $s$-term approximation error are as follows:

\begin{defi}
	The weighted $l_p$-space is defined as $$l_{\boldsymbol{\omega},p}:=\left\{\mathbf{x}=(x_j)_{j\in \Lambda},||\mathbf{x}||_{\boldsymbol{\omega},p}:=\left(\mathop{\sum}\limits_{j\in\Lambda}\omega_j^{2-p}|x_j|^p\right)^{\frac{1}{p}}<\infty\right\},\quad 0<p\leqslant 2,$$ 
    from which we also have $||\mathbf{x}||_{\boldsymbol{\omega},0}=\mathop{\sum}_{\{j:x_j\neq 0\}}\omega_j^2$ and the weighted cardinality of any subset $S\subset\Lambda$ is given by $\omega(S):=\mathop{\sum}_{j\in S}\omega_j^2$. In this space, the weighted best $s$-term approximation error of a vector $\mathbf{x}$ under the $l_{\boldsymbol{\omega},p}$-norm is $\sigma_s(\mathbf{x})_{\boldsymbol{\omega},p}=\mathop{\inf}\limits_{\{\mathbf{z}:||\mathbf{z}||_{\boldsymbol{\omega},0}\leqslant s\}}||\mathbf{x}-\mathbf{z}||_{\boldsymbol{\omega},p}$.
\end{defi}

It is clear that by assigning a distinct weight to each component of a vector, we arrive at the concept of a weighted $l_p$-space, along with the extension of $l_0$-norm and best $s$-term approximation error within the $l_p$-space. The weighted $l_1$-norm of a vector $\mathbf{z}\in\mathbb{R}^N$ is particularly frequently used and is given by $||\mathbf{z}||_{\boldsymbol{\omega},1}=\mathop{\sum}_{i=1}^N\omega_i|z_i|$.

A compressive continuous function is characterized as a continuous function for which the best $s$-term approximation error, or the weighted best $s$-term approximation error of the coefficient vector in its generalized Fourier expansion with respect to an orthonormal polynomial basis, decreases rapidly as $s$ increases.

For polynomial function $g(x)=\mathop{\sum}_{k=1}^Nc_k\psi_k(x)$, its quasi-norm $||\cdot||_{A_q}$ is defined as
$$||g||_{A_q}:=||\mathbf{c}||_q,\quad q>0.$$

Additionally, suppose that the sampling points $x_1,x_2,\dots,x_m$ are chosen randomly in $I$ according to the probability density function given by the probability measure $\nu$. Then, using the $m$ sampled values of the sparse polynomial $g(x)$, we can recover the $s$-sparse polynomial $g(x)$ through an $l_1$-minimization problem. The subsequent conclusion pertain to the recovery guarantee and approximation error associated with the recovery of sparse polynomial $g(x)$ based on $l_1$-minimization problem:

\begin{thm}[\cite{HR2012}]\label{three}
	Assume that $N,m,s\in\mathbb{N}$ are given such that
    \begin{equation}\label{mboundbegin}
	m\geqslant CK^2s\log^3(s)\log(N).
    \end{equation}
Suppose that $m$ sampling points $\{x_1,x_2,\dots,x_m\}$ are independently and randomly selected according to the probability density function given by the probability measure $\nu$, and $\{\psi_k\}_{k=1}^N$ be a set of polynomials that are uniformly bounded and orthonormal with respect to the measure function $\nu$, with a uniform bound $K$. The sampling matrix is defined as $A=(A_{j,k})$ where $A_{j,k}=\psi_{k}(x_j)$ for $j=1,\dots,m$ and $k=1,\dots,N$. Then, with probability at least $1-N^{-\gamma \log^3(s)}$, the following holds: for any polynomial $g(x)=\mathop{\sum}_{k=1}^Nc_k\psi_k(x)$, if noisy sampled values $\mathbf{y}=\left(g(x_1)+\eta_1,\dots,g(x_m)+\eta_m\right)^T=A\mathbf{c}+\boldsymbol{\eta}$ are obtained with $||\boldsymbol{\eta}||_\infty\leqslant\varepsilon$, the coefficient vector $\mathbf{c}=(c_1,\dots,c_N)^T$ can be recovered by solving the $l_1$-minimization problem with inequality constraint:
\begin{equation}\label{l1minibegin}
\mathbf{c}^\#=\arg\mathop{\min}\limits_{\mathbf{z}\in\mathbb{R}^N}||\mathbf{z}||_1\quad \mathrm{subject}\hskip 0.5em \mathrm{to}\hskip 0.5em||A\mathbf{z}-\mathbf{y}||_2\leqslant\sqrt{m}\varepsilon,
\end{equation}
Precisely,
	$$\big|\big|\mathbf{c}-\mathbf{c}^{\#}\big|\big|_2\leqslant C_1\frac{\sigma_s(\mathbf{c})_1}{\sqrt{s}}+C_2\varepsilon,$$
	as well as
	$$\big|\big|\mathbf{c}-\mathbf{c}^{\#}\big|\big|_1\leqslant D_1\sigma_s(\mathbf{c})_1+D_2\sqrt{s}\varepsilon,$$
	where $C,C_1,C_2,D_1,D_2,\gamma>0$ are universal constants.
\end{thm}

Similarly, drawing from the conclusions about the weighted Restricted Isometry Property (RIP) of the measurement matrix $A$ and the error bounds for the reconstructed vector obtained by solving weighted $l_1$-minimization problem as detailed in \cite{HR2016}, we can infer the following conclusion concerning the recovery guarantee and approximation error associated with the recovery of sparse polynomial $g(x)$ when using weighted $l_1$-minimization approach:

\begin{thm}[\cite{HR2016}]\label{weight}
	For the orthonormal system $\{\psi_j\}_{j=1}^{N}$, let the weight vector $\boldsymbol{\omega}$ have components $\omega_j\geqslant||\psi_j||_\infty$, for $j=1,\dots,N$, and let $s\geqslant 2||\boldsymbol{\omega}||_\infty^2$. By selecting the sampling points $\{x_1,x_2,\dots,x_m\}$ independently and randomly selected according to the probability measure $\nu$ as the probability density function, we construct the matrix $A\in \mathbb{R}^{m\times N}$ such that $A_{j,k}=\psi_{k}(x_j)$, for $j=1,\dots,m$ and $k=1,\dots,N$. Then, for any function $g(x)=\mathop{\sum}_{i=1}^{N}c_i\psi_i(x)$ generated by this orthonormal system in the linear space, denoting the noisy sampled values at the sampling locations as $\mathbf{y}=\left(g(x_1)+\eta_1,\dots,g(x_m)+\eta_m\right)^T=A\mathbf{c}+\boldsymbol{\eta}$ with $||\boldsymbol{\eta}||_\infty\leqslant\varepsilon$, if the number of samples satisfies
	\begin{equation}\label{mbound43}
		m\geqslant Cs\log^3(s)\log(N),
	\end{equation}	
	then with a probability of at least $1-N^{-\log^3(3s)}$, solving the weighted $l_1$-minimization problem
	\begin{equation}\label{BOSq}
\mathop{\min}\limits_{\mathbf{z}\in\mathbb{R}^N}||\mathbf{z}||_{\boldsymbol{\omega},1}\quad \mathrm{subject}\hskip 0.5em \mathrm{to}\hskip 0.5em||A\mathbf{z}-\mathbf{y}||_2\leqslant\sqrt{m}\varepsilon    
	\end{equation}
	yields a solution $\mathbf{c}^\#$ that satisfies the error estimates
	$$
	\left\{
	\begin{aligned}
		&\big|\big|\mathbf{c}-\mathbf{c}^\#\big|\big|_{\boldsymbol{\omega},1}\leqslant F_1\sigma_s(\mathbf{c})_{\boldsymbol{\omega},1}+G_1\sqrt{s}\varepsilon,\\
		&\big|\big|\mathbf{c}-\mathbf{c}^\#\big|\big|_2\leqslant F_2\frac{\sigma_s(\mathbf{c})_{\boldsymbol{\omega},1}}{\sqrt{s}}+G_2\varepsilon,
	\end{aligned}
	\right.
	$$
	where $C,F_1,G_1,F_2,G_2>0$ are constants.
\end{thm}

\section{The $L_\infty$-norm Estimate of the Approximation Error of Sparse Polynomial to Continuous Function}

In this section, we employ random sampling of the continuous function $f(x)$ to construct a sparse polynomial approximation $p^\#(x)$ of $f(x)$ \textit{via} an $l_1$-minimization problem and estimate the magnitude of the approximation error $||f(x)-p^\#(x)||_\infty$.

\begin{thm}\label{thm4}
Suppose $0<q<1$, and let $N,m,s \in\mathbb{N}$ satisfy
\begin{equation}\label{mmm}
    m\geqslant CK^2s\log^3(s)\log(N+1),
\end{equation}
$\{\psi_0,\psi_1,\dots,\psi_{N}\}$ is a set of uniformly bounded orthonormal polynomials on interval $I=[a,b]$ with respect to the probability measure $\nu(x)$, with a uniform bound $K$. The degrees of these polynomials strictly increase and do not exceed $N$. The $m$ sampling points $x_1,\dots,x_m$ are randomly selected on $I$ according to the probability density function given by the probability measure $\nu(x)$ and satisfy independent and identical distribution. The matrix $A$ is defined with entries $A_{j,k}=\psi_{k}(x_j)$ for $j=1,...,m$ and $k=0,1,...,N$. Then, for any continuous function $f(x)\in C[a,b]$, let the vector of sampled values be $\mathbf{y}=(f(x_1),\dots,f(x_m))^T$. With a probability of at least $1-(N+1)^{-\gamma \log^3(s)}$, the following holds: if $\mathbf{c}^{\#}=(c_0^\#,c_1^\#,\dots,c_{N}^\#)^T$ is the solution to the problem
	\begin{equation}\label{Infq}
	\mathop{\min}\limits_{\mathbf{z}\in \mathbb{R}^{N+1}}||\mathbf{z}||_1\quad\mathrm{subject}\hskip 0.5em  \mathrm{to}\hskip 0.5em  ||A\mathbf{z}-\mathbf{y}||_2\leqslant\Theta\sqrt{m}E_{N},
	\end{equation}
	where $\Theta$ is a predetermined positive integer, and let $p^\#(x)=\mathop{\sum}_{k=0}^{N}c_k^\#\psi_k(x)$, then the error estimate is given by
	\begin{equation}\label{errorInf}
		||f-p^\#||_\infty\leqslant K\cdot D_1\sigma_s(\mathbf{c_{opt}})_1+\left(K\cdot D_2\Theta\sqrt{s}+1\right)E_{N}
	\end{equation}
	or
	\begin{equation}\label{errorInf2}
		||f-p^\#||_\infty\leqslant K\cdot D_1||p_{opt}||_{A_q}s^{1-\frac{1}{q}}+\left(K\cdot D_2\Theta\sqrt{s}+1\right)E_{N},
	\end{equation}
	where $C,D_1,D_2,\gamma>0$ are constants in the theorem. $p_{opt}(x)$ is the best uniform approximation polynomial of $f(x)$ in polynomial space $\Pi_N(I)=\mathrm{span}\{\psi_0,\psi_1,\dots,\psi_{N}\}$, $\mathbf{c_{opt}}$ is the coefficient vector of $p_{opt}$ in orthonormal basis $\{\psi_j\}_{j=0}^N$, and $E_{N}$ is the best uniform approximation error of $f(x)$ in $\Pi_{N}(I)$.
\end{thm}

{\it Proof.} Since $p_{opt}$ is the best uniform approximation polynomial in the finite-dimensional polynomial space $\Pi_N(I)=\mathrm{span}\{\psi_0,\psi_1,\dots,\psi_{N}\}$ for $f(x)$, it follows that $p_{opt}=\mathop{\sum}_{k=0}^Nc_{opt,k}\psi_k(x)$ and satisfies $||f-p_{opt}||_\infty\leqslant E_N=\mathop{\inf}\limits_{p\in\Pi_{N}(I)}\{\Delta(p)\}$. Considering $f(x_l)$ as a noisy sample of $p_{opt}(x_l)$, we have
	$$f(x_l)=p_{opt}(x_l)+\eta_l,\quad l=1,\dots,m,$$
	from which it follows that
	$$|\eta_l|=|f(x_l)-p_{opt}(x_l)|\leqslant E_N\leqslant \Theta E_N,\quad l=1,\dots,m,$$
	where $\Theta$ is any positive integer. Furthermore,
	$$||A\mathbf{c_{opt}}-\mathbf{y}||_2\leqslant\Theta\sqrt{m}E_N,\quad \mathbf{y}=\left(f(x_1),\dots,f(x_m)\right)^T.$$
	
	Let $\mathbf{c}^\#=(c_0^\#,c_1^\#,\dots,c_{N}^\#)^T$ be the solution to the optimization problem (\ref{Infq}). According to Theorem \ref{three}, when $m$ satisfies equation (\ref{mmm}), we have
	\begin{equation}\label{Infl1norm}
		||\mathbf{c_{opt}}-\mathbf{c}^\#||_1\leqslant D_1\sigma_s(\mathbf{c_{opt}})_1+D_2\Theta\sqrt{s}E_{N}.
	\end{equation}
    
	Using $\mathbf{c}^\#$ to generate the sparse polynomial $p^\#(x)=\mathop{\sum}_{k=0}^{N}c_k^\#\psi_k(x)$, then
	\begin{equation}\label{Inffp}
		\begin{split}
			||f-p^\#||_\infty&\leqslant ||f-p_{opt}||_\infty+||p_{opt}-p^{\#}||_\infty\\
			&=E_{N}+\bigg|\bigg|\mathop{\sum}\limits_{k=0}^{N}(c_k^\#-c_{opt,k})\psi_k(x)\bigg|\bigg|_\infty\\
			&\leqslant E_{N}+\mathop{\sum}\limits_{k=0}^{N}\big|c_k^\#-c_{opt,k}\big|\cdot ||\psi_k(x)||_\infty\\
			&\leqslant E_{N}+
			K\cdot||\mathbf{c}^\#-\mathbf{c_{opt}}||_1.
		\end{split}
	\end{equation}
	Substituting equation (\ref{Infl1norm}) into equation (\ref{Inffp}), we get:
	\begin{equation}\label{Inffp2}
		\begin{split}
			||f-p^\#||_\infty&\leqslant E_{N}+K\left(D_1\sigma_s(\mathbf{c_{opt}})_1+D_2\Theta\sqrt{s}E_{N}\right)\\
			&=K\cdot D_1\sigma_s(\mathbf{c_{opt}})_1+\left(K\cdot D_2\Theta\sqrt{s}+1\right)E_{N}.
		\end{split}
	\end{equation}
	According to the Stechkin estimate,
	\begin{equation}\label{Stechkin}
    \sigma_s(\mathbf{c_{opt}})_1\leqslant s^{1-\frac{1}{q}}||\mathbf{c_{opt}}||_q,\quad 0<q<1,
	\end{equation}
	substituting (\ref{Stechkin}) into (\ref{Inffp2}), and using the definition of $||p||_{A_q}$, we obtain
	$$||f-p^\#||_\infty\leqslant K\cdot D_1||p_{opt}||_{A_q}s^{1-\frac{1}{q}}+\left(K\cdot D_2\Theta\sqrt{s}+1\right)E_{N},\quad 0<q<1.$$
	Thus, Theorem \ref{thm4} is proven.$\hfill \Box$ \@
	
\begin{rem}
	If $\mathbf{c_{opt}}$ is compressible or $||p_{opt}||_{A_q}\leqslant 1$, then by Theorem \ref{thm4}, as $s$ increases, $p^\#(x)$ gradually approximates $f(x)$. If $f(x)\in\Pi_N(I)$ and is an $s$-sparse polynomial, we have $E_N=\sigma_s(\mathbf{c_{opt}})_1=0$, leading to $p^\#(x)=f(x)$, and thus $f(x)$ is exactly reconstructed.
\end{rem}

Using Corollary \ref{coro1}, we can infer the subsequent Corollary \ref{coro3} from Theorem \ref{thm4}.

\begin{deduction}\label{coro3}
	If $f(x)\in Lip_M^\alpha(0<\alpha\leqslant 1)$, then
	$$||f-p^\#||_\infty\leqslant K\cdot D_1||p_{opt}||_{A_q}s^{1-\frac{1}{q}}+\left(K\cdot D_2\Theta+1\right)CM\frac{\sqrt{s}}{N^\alpha},\quad 0<q<1,$$
	where $C=12\left(\frac{b-a}{2}\right)^\alpha.$	
\end{deduction}

\begin{rem}
	Setting $\sqrt{s}/{N^\alpha}=s^{1-1/q}$, we derive that $N=s^{(1/q-1/2)/\alpha}$. By choosing $N=\lceil s^{(1/q-1/2)/\alpha}\rceil$, we obtain
	\begin{equation}\label{38}
	||f-p^\#||_\infty\leqslant \left(K\cdot D_1||p_{opt}||_{A_q}+\left(K\cdot D_2\Theta+1\right)CM\right)s^{1-\frac{1}{q}},\quad 0<q<1.
	\end{equation}
\end{rem}

Using Corollary \ref{coro2}, we can infer the subsequent Corollary \ref{coro4} from Theorem \ref{thm4}.

\begin{deduction}\label{coro4}
	If function $f(x)$ has a bounded derivative $f'(x)$ such that $|f'(x)|\leqslant M_{1}$, then
	$$||f-p^\#||_\infty\leqslant K\cdot D_1||p_{opt}||_{A_q}s^{1-\frac{1}{q}}+6(b-a)M_{1}\left(K\cdot D_2\Theta+1\right)\frac{\sqrt{s}}{N},\quad 0<q<1.$$	
\end{deduction}

\begin{rem}
	Setting $\sqrt{s}/{N}=s^{1-1/q}$, we derive that $N=s^{1/q-1/2}$. By choosing $N=\lceil s^{1/q-1/2}\rceil$, we obtain
	\begin{equation}\label{39}
	||f-p^\#||_\infty\leqslant \left(K\cdot D_1||p_{opt}||_{A_q}+6(b-a)M_{1}\left(K\cdot D_2\Theta+1\right)\right)s^{1-\frac{1}{q}},\quad 0<q<1.
	\end{equation}
\end{rem}

By applying Theorem \ref{thm2} and Theorem \ref{thm4}, we arrive at the subsequent Theorem \ref{thm5}.

\begin{thm}\label{thm5}
	If $f(x)$ has a bounded $p$-th derivative $f^{(p)}(x)$ such that $|f^{(p)}(x)|\leqslant M_p$, then
	$$||f-p^\#||_\infty\leqslant K\cdot D_1||p_{opt}||_{A_q}s^{1-\frac{1}{q}}+C_p(b-a)^pM_p\left(K\cdot D_2\Theta+1\right)\frac{\sqrt{s}}{N^p},\quad 0<q<1,$$
	where the constant $C_p$ depends only on $p\geqslant 2$.	
\end{thm} 

\begin{rem}
	Setting $\sqrt{s}/{N^p}=s^{1-1/q}$, we derive that $N=s^{(1/q-1/2)/p}$. By choosing $N=\lceil s^{(1/q-1/2)/p}\rceil$, we obtain
	\begin{equation}\label{310}
		||f-p^\#||_\infty\leqslant \left(K\cdot D_1||p_{opt}||_{A_q}+C_p(b-a)^pM_p\left(K\cdot D_2\Theta+1\right)\right)s^{1-\frac{1}{q}},\quad 0<q<1.
	\end{equation}
\end{rem}

From equations (\ref{38}) to (\ref{310}), we can see that in order to achieve the same upper bound on error, if the smoothness order $p$ of $f(x)$ is higher, the required degree $N$ of the sparse polynomial approximation is smaller. This suggests that the greater the smoothness of $f(x)$, the more readily it can be approximated by a sparse polynomial.

\section{The $L_2$-norm Estimate of the Approximation Error of Sparse Polynomial to Continuous Function}
Let $p_{opt}$ represent the best square approximation polynomial of the continuous function $f(x)$ in polynomial space $\Pi_N(I)=\mathrm{span}\{\psi_0,\dots,\psi_N\}$, such that, $||f-p_{opt}||_2=\mathop{\inf}\limits_{p\in\Pi_N(I)}\{\Delta(p)\}$, where $\Delta(p)=||f-p||_2.$

In this section, we persist in employing random sampling of the continuous function $f(x)$ to construct a sparse approximation polynomial $p^\#(x)$ through $l_1$-minimization problem or weighted $l_1$-minimization problem, and estimate the magnitude of the approximation error $||f(x)-p^\#(x)||_2$. Initially, we continue with the method from the previous section that utilizes an $l_1$-minimization problem to generate a sparse polynomial, which leads us to the following theorem:

\begin{thm}\label{l1Square}
	Suppose $0<q<1$, and let $N,m,s \in\mathbb{N}$ satisfy
	\begin{equation}\label{mmmm}
		m\geqslant CK^2s\log^3(s)\log(N+1).
	\end{equation}
	Let $\{\psi_0,\psi_1,\dots,\psi_{N}\}$ be a set of uniformly bounded orthonormal polynomials defined on $I=[a,b]$ with respect to the probability measure $\nu(x)$, with a uniform bound $K$. Their degrees strictly increase and do not exceed $N$. The $m$ sampling points $x_1,\dots,x_m$ are randomly selected on $I$ according to the probability density function given by $\nu(x)$ and are independently and identically distributed. The matrix $A$ is defined with elements $A_{j,k}=\psi_{k}(x_j)$ for $j=1,...,m$ and $k=0,1,...,N$. Then, for any continuous function $f(x)\in C[a,b]$ on $I$, let the sampled value vector be $\mathbf{y}=(f(x_1),\dots,f(x_m))^T$. With a probability of at least $1-(N+1)^{-\gamma \log^3(s)}$, the following holds: if $\mathbf{c}^{\#}=(c_0^\#,c_1^\#,\dots,c_{N}^\#)^T$ is the solution to the problem
	\begin{equation}\label{Squareq}
	\mathop{\min}\limits_{\mathbf{z}\in \mathbb{R}^{N+1}}||\mathbf{z}||_1\quad\mathrm{subject}\hskip 0.5em  \mathrm{to}\hskip 0.5em  ||A\mathbf{z}-\mathbf{y}||_2\leqslant\Theta\sqrt{m}||T_{N}||_\infty,    
	\end{equation}
where $\Theta$ is a predetermined positive integer, and let $p^\#=\mathop{\sum}_{k=0}^{N}c_k^\#\psi_k$, then the error estimate is given by
	\begin{equation}\label{errorSquare}
		||f-p^\#||_2\leqslant ||T_{N}||_2+C_1\frac{\sigma_s(\mathbf{c_{opt}})_1}{\sqrt{s}}+C_2\Theta ||T_{N}||_\infty,
	\end{equation}
	or
	\begin{equation}\label{errorSquare2}
		||f-p^\#||_2\leqslant ||T_{N}||_2+C_1||f||_{A_q}s^{\frac{1}{2}-\frac{1}{q}}+C_2\Theta ||T_{N}||_\infty,
	\end{equation}
	where $C,C_1,C_2,\gamma>0$ are constants in the theorem. $\mathbf{c_{opt}}$ is the coefficient vector of the best square approximation polynomial $p_{opt}$ of $f$ in $\Pi_N(I)=\mathrm{span}\{\psi_0,\psi_1,\dots,\psi_N\}$, and $T_{N}(x)=f(x)-p_{opt}(x)$.
\end{thm}

{\it Proof.} Let $T_{N}(x)=f(x)-p_{opt}(x)$, and denote $\mathbf{c_{opt}}=(c_{opt,0},\dots,c_{opt,N})^T$, so that $p_{opt}=\mathop{\sum}_{k=0}^{N}c_{opt,k}\psi_{k}(x)$. Since $p_{opt}$ is the best square approximation polynomial of $f(x)$ in $\Pi_N(I)$, and $f(x)-p_{opt}(x)=T_{N}(x)$, we know that if we consider the sampled values $f(x_i)$ as noisy samples of $p_{opt}(x_i)$ for $i=1,\dots,m$：
\begin{equation}\label{SquareNoise}
	f(x_i)=p_{opt}(x_i)+\eta_i,\quad i=1,\dots,m.
\end{equation}
Then $|\eta_i|=|f(x_i)-p_{opt}(x_i)|\leqslant ||T_{N}||_\infty,i=1,\dots,m$, meaning the sampling noise vector $\boldsymbol{\eta}=(\eta_1,\dots,\eta_m)^T$ has an upper bound $||\boldsymbol{\eta}||_\infty\leqslant \Theta ||T_{N}||_\infty$, where $\Theta$ is a positive integer. According to the construction of the measurement matrix in the theorem, $\mathbf{y}=A\mathbf{c_{opt}}+\boldsymbol{\eta}$, and thus, by Theorem \ref{three}, the solution $\mathbf{c}^\#=(c_0^\#,c_1^\#,\dots,c_{N}^\#)^T$ obtained from solving the optimization problem (\ref{Squareq}) satisfies the following error estimate with a probability of at least $1-(N+1)^{-\gamma \log^3(s)}$:
\begin{equation}\label{Squarel1norm}
	||\mathbf{c}^\#-\mathbf{c_{opt}}||_2\leqslant C_1\frac{\sigma_s(\mathbf{c_{opt}})_1}{\sqrt{s}}+C_2\Theta ||T_{N}||_\infty,
\end{equation}
Let $p^{\#}(x)=\mathop{\sum}_{k=0}^{N}c_k^\#\psi_k$, and consider the approximation error $||f-p^{\#}||_2$:
\begin{equation}\label{Squarefp}
	\begin{split}
		||f-p^\#||_2&\leqslant ||f-p_{opt}||_2+||p_{opt}-p^{\#}||_2\\
		&=||T_{N}||_2+\bigg|\bigg|\mathop{\sum}\limits_{k=0}^{N}(c_k^\#-c_{opt,k})\psi_k(x)\bigg|\bigg|_2\\
		&\leqslant ||T_{N}||_2+\left(\int_a^b\left(\mathop{\sum}\limits_{k=0}^{N}(c_k^\#-c_{opt,k})\psi_k(x)\right
		)^2\rho(x)dx\right)^{\frac{1}{2}},
	\end{split}
\end{equation}
By the orthonormality of $\{\psi_k\}_{k=0}^{N}$, equation (\ref{Squarefp}) becomes:
\begin{equation}\label{Squarefp2}
	\begin{split}
		||f-p^\#||_2&\leqslant ||T_{N}||_2+\left(\mathop{\sum}\limits_{k=0}^{N}(c_k^\#-c_{opt,k})^2\int_a^b\psi_k^2(x)\rho(x)dx\right)^{\frac{1}{2}}\\
		&=||T_{N}||_2+\left(\mathop{\sum}\limits_{k=0}^{N}(c_k^\#-c_{opt,k})^2\right)^{\frac{1}{2}}\\
		&=||T_{N}||_2+||\mathbf{c}^\#-\mathbf{c_{opt}}||_2,
	\end{split}
\end{equation}   
Substituting the result of (\ref{Squarel1norm}) into (\ref{Squarefp2}), we get:
$$
||f-p^\#||_2\leqslant ||T_{N}||_2+C_1\frac{\sigma_s(\mathbf{c_{opt}})_1}{\sqrt{s}}+C_2\Theta ||T_{N}||_\infty,
$$
Thus, the error estimate (\ref{errorSquare}) in the theorem is proven. Furthermore, using the Stechkin inequality and $||p_{opt}||_{A_q}\leqslant ||f||_{A_q}$, the above inequality can be extended to
$$
||f-p^\#||_2\leqslant ||T_{N}||_2+C_1||f||_{A_q}s^{\frac{1}{2}-\frac{1}{q}}+C_2\Theta ||T_{N}||_\infty,
$$
which completes the proof of Theorem \ref{l1Square}.$\hfill \Box$ \@

\begin{rem}
	If $\mathbf{c_{opt}}$ is compressible or $||p_{opt}||_{A_q}\leqslant 1$, then by Theorem \ref{l1Square}, as $s$ increases, $p^\#(x)$ gradually approximates $f(x)$. When $f(x)\in\Pi_N(I)$ and is an $s$-sparse polynomial, we have $||T_N||_2=||T_N||_\infty=\sigma_s(\mathbf{c_{opt}})_1=0$, leading to $p^\#(x)=f(x)$, and thus $f(x)$ is exactly reconstructed.
\end{rem}

Theorem \ref{thm4} and Theorem \ref{l1Square} both utilize $l_1$-minimization problem to construct sparse polynomial approximations for continuous functions. Moving forward, we will explore constructing analogous theorems based on weighted $l_1$-minimization problem.

Leveraging the conclusions from Theorem \ref{weight} and adhering to the proof methodology of Theorem \ref{thm4}, we can deduce a method for sparse polynomial approximation of functions that incorporates the best uniform approximation error and the weighted $l_1$-minimization approach:

\begin{thm}\label{l1Infweight}
	Suppose $N,m,s \in\mathbb{N}$ satisfy
	\begin{equation}
		m\geqslant Cs\log^3(s)\log(N+1).
	\end{equation}
	Let $\{\psi_0,\psi_1,\dots,\psi_{N}\}$ be a set of orthonormal polynomials defined on $I=[a,b]$ with respect to the probability measure $\nu(x)$, with their degrees strictly increasing and not exceeding $N$. Suppose the selected weights are $\boldsymbol{\omega}=(\omega_0,\omega_1,\dots,\omega_N)^T$ where each component satisfies $\omega_j\geqslant||\psi_j||_\infty$, for $j=0,1,\dots,N$, and $s\geqslant 2||\boldsymbol{\omega}||_\infty^2$. The $m$ sampling points $x_1,\dots,x_m$ are randomly selected on $I$ according to the probability density function given by $\nu(x)$ and are independently and identically distributed. The matrix $A$ is defined with elements $A_{j,k}=\psi_{k}(x_j)$ for $j=1,...,m$ and $k=0,1,...,N$. Then, for any continuous function $f(x)\in C[a,b]$ on $I$, let the sampled value vector be $\mathbf{y}=(f(x_1),\dots,f(x_m))^T$. With a probability of at least $1-(N+1)^{-\log^3(3s)}$, the following holds: if $\mathbf{c}^{\#}=(c_0^\#,c_1^\#,\dots,c_{N}^\#)^T$ is the solution to the problem
	\begin{equation}\label{Infqweight}
	\mathop{\min}\limits_{\mathbf{z}\in \mathbb{R}^{N+1}}||\mathbf{z}||_{\boldsymbol{\omega},1}\quad\mathrm{subject}\hskip 0.5em  \mathrm{to}\hskip 0.5em  ||A\mathbf{z}-\mathbf{y}||_2\leqslant\Theta\sqrt{m}E_{N},
	\end{equation}
	where $\Theta$ is a predetermined positive integer, and let $p^\#(x)=\mathop{\sum}_{k=0}^{N}c_k^\#\psi_k$, then the error estimate is given by
	\begin{equation}\label{errorInfweight}
		||f-p^\#||_2\leqslant ||f-p_{opt}||_2+ F_2\frac{\sigma_s(\mathbf{c_{opt}})_{\boldsymbol{\omega},1}}{\sqrt{s}}+\Theta G_2E_{N},
	\end{equation}
	where $C,F_2,G_2>0$ are constants in the theorem. $\mathbf{c_{opt}}$ is the coefficient vector of the best uniform approximation polynomial $p_{opt}$ of $f$ in $\Pi_N(I)=\mathrm{span}\{\psi_0,\psi_1,\dots,\psi_N\}$, and $E_{N}$ is the best uniform approximation error of $f$ in $\Pi_{N}(I)$.
\end{thm}

\begin{rem}
	The error estimate (\ref{errorInfweight}) is derived as follows:
	\begin{equation}\label{Squarefpweight}
		\begin{split}
			||f-p^\#||_2&\leqslant ||f-p_{opt}||_2+||p_{opt}-p^{\#}||_2\\
			&=||f-p_{opt}||_2+\bigg|\bigg|\mathop{\sum}\limits_{k=0}^{N}(c_k^\#-c_{opt,k})\psi_k(x)\bigg|\bigg|_2\\
			&\leqslant ||f-p_{opt}||_2+\left(\int_a^b\left(\mathop{\sum}\limits_{k=0}^{N}(c_k^\#-c_{opt,k})\psi_k(x)\right
			)^2\rho(x)dx\right)^{\frac{1}{2}}.
	\end{split}
    \end{equation}
    By the orthonormality of $\{\psi_k\}_{k=0}^{N}$, equation (\ref{Squarefpweight}) becomes:
\begin{equation}\label{Squarefpweight2}
		\begin{split}
			||f-p^\#||_2
    &\leqslant ||f-p_{opt}||_2+\left(\mathop{\sum}\limits_{k=0}^{N}(c_k^\#-c_{opt,k})^2\right)^{\frac{1}{2}}\\
			&=||f-p_{opt}||_2+||\mathbf{c}^\#-\mathbf{c_{opt}}||_2\leqslant||f-p_{opt}||_2+ F_2\frac{\sigma_s(\mathbf{c_{opt}})_{\boldsymbol{\omega},1}}{\sqrt{s}}+\Theta G_2E_{N}.
		\end{split}
        \end{equation}
\end{rem}

Comparing Theorem \ref{l1Infweight} with Theorem \ref{thm4}, we notice that the former does not require the orthonormal polynomial system to be uniformly bounded, and the uniform bound does not appear in the lower bound for the number of samples needed. The derived error estimate is for $L_2$-norm, but the term $||f-p_{opt}||_2$ remains very small if $N$ is sufficiently large. Together with the already small $E_N$ and the fact that $\frac{\sigma_s(\mathbf{c_{opt}})_{\boldsymbol{\omega},1}}{\sqrt{s}}$ is not significant when $s$ is large, the overall upper bound on the error is also very small. On this basis, if $\mathbf{c_{opt}}$ is compressible, then according to Theorem \ref{l1Infweight}, as $s$ increases, $p^\#(x)$ gradually approximates $f(x)$. If $f(x)\in\Pi_N(I)$ and is an $s$-sparse polynomial, then $E_N=0$ and $f-p_{opt}=0$, and $\sigma_s(\mathbf{c_{opt}})_{\boldsymbol{\omega},1}$ is very close to $0$, thus the error between $p^\#(x)$ and $f(x)$ is very small, and $f(x)$ is approximately exactly reconstructed.

Following a similar proof process as Theorem \ref{l1Square} and combining it with Theorem \ref{weight}, we arrive at the following theorem:

\begin{thm}\label{l1Squareweight}
	Suppose $N,m,s \in\mathbb{N}$ satisfy
	\begin{equation}
		m\geqslant Cs\log^3(s)\log(N+1)
	\end{equation}
	Let $\{\psi_0,\psi_1,\dots,\psi_{N}\}$ be a set of orthonormal polynomials defined on $I=[a,b]$ with respect
to the probability measure $\nu(x)$, with their degrees strictly increasing and not exceeding $N$. Suppose the selected weights are $\boldsymbol{\omega}=(\omega_0,\omega_1,\dots,\omega_N)^T$ where each component satisfies $\omega_j\geqslant||\psi_j||_\infty$, for $j=0,1,\dots,N$, and $s\geqslant 2||\boldsymbol{\omega}||_\infty^2$. The $m$ sampling points $x_1,\dots,x_m$ are randomly selected on $I$ according to the probability density function given by $\nu(x)$ and are independently and identically distributed. The matrix $A$ is defined with elements $A_{j,k}=\psi_{k}(x_j)$ for $j=1,...,m$ and $k=0,1,...,N$. Then, for any continuous function $f(x)\in C[a,b]$ on $I$, let the sampled value vector be $\mathbf{y}=(f(x_1),\dots,f(x_m))^T$. With a probability of at least $1-(N+1)^{-\log^3(3s)}$, the following holds: if $\mathbf{c}^{\#}=(c_0^\#,c_1^\#,\dots,c_{N}^\#)^T$ is the solution to the problem
	\begin{equation}\label{Squareqweight}
	\mathop{\min}\limits_{\mathbf{z}\in \mathbb{R}^{N+1}}||\mathbf{z}||_{\boldsymbol{\omega},1}\quad\mathrm{subject}\hskip 0.5em  \mathrm{to}\hskip 0.5em  ||A\mathbf{z}-\mathbf{y}||_2\leqslant\Theta\sqrt{m}||T_{N}||_\infty,    
	\end{equation}
	where $\Theta$ is a predetermined positive integer, and let $p^\#=\mathop{\sum}_{k=0}^{N}c_k^\#\psi_k$, then the error estimate is given by
\begin{equation}\label{errorSquareweight}
		||f-p^\#||_2\leqslant ||T_{N}||_2+F_2\frac{\sigma_s(\mathbf{c_{opt}})_{\boldsymbol{\omega},1}}{\sqrt{s}}+\Theta G_2||T_{N}||_\infty,
	\end{equation}
	where $C,F_2,G_2>0$ are constants in the theorem. $\mathbf{c_{opt}}$ is the coefficient vector of the best square approximation polynomial $p_{opt}$ of $f$ in $\Pi_N(I)=\mathrm{span}\{\psi_0,\psi_1,\dots,\psi_N\}$, and $T_{N}$ represents the best square approximation error function $f-p_{opt}$ of $f$ within the space $\Pi_{N}(I)$.
\end{thm}

Comparing Theorem \ref{l1Squareweight} with Theorem \ref{l1Square}, we notice that the former does not require the orthonormal polynomial system to be uniformly bounded, and the error estimates derived from both theorems are similar. If $\mathbf{c_{opt}}$ is compressible, by Theorem \ref{l1Squareweight}, as $s$ increases, $p^\#(x)$ gradually approximates $f(x)$. When $f(x)\in\Pi_N(I)$ and is an $s$-sparse polynomial, we have $||T_N||_2=||T_N||_\infty=0$ and $\sigma_s(\mathbf{c_{opt}})_{\boldsymbol{\omega},1}$ is negligible, which implies that $p^\#(x)$ closely approximates $f(x)$, leading to an approximate exact reconstruction of $f(x)$.

\section{Numerical Experiments}
In this section, we conduct numerical experiments based on the theories previously discussed. We utilize the YALL1 algorithm package on MATLAB, designed by Yang and Zhang, to solve the $l_1$-minimization problem and the weighted $l_1$-minimization problem \cite{JY2011}. For these experiments, we employ the normalized first-kind Chebyshev polynomial system on the interval $[-1,1]$ as our basis. The normalization process is detailed as follows: For the first-kind Chebyshev orthogonal polynomials $\{T_n(x)\}_{n=0}^\infty$ on $[-1,1]$, defined by $T_n(x)=\cos(n\cdot \arccos(x))$, we keep $T_0(x)$ unaltered, and multiply each subsequent term $T_n(x)$ by $\sqrt{2}$ for $n=1,2,\dots$, yielding a new set of orthogonal polynomials $\{C_n(x)\}_{n=0}^\infty$. Given the property of the first-kind Chebyshev orthogonal polynomials:
\begin{equation}\label{Cheby}
	\int_{-1}^{1}\frac{T_j(x)T_k(x)}{\sqrt{1-x^2}}dx=\left\{
	\begin{array}{l}
		\pi,\quad j=k=0,\\
		\frac{\pi}{2},\quad j=k\neq 0,\\
		0,\quad j\neq k,
	\end{array}
	\right.
\end{equation}
we can infer that $\{C_n(x)\}_{n=0}^\infty$ forms an orthonormal polynomial system with respect to the probability measure $\nu_C(x)=\frac{1}{\pi\sqrt{1-x^2}}$, which corresponds to the probability density function of the arcsine distribution on $[-1,1]$. This system is referred to as the normalized first-kind Chebyshev polynomial system and satisfies:
\begin{equation}
	\int_{-1}^1\frac{C_j(x)C_k(x)}{\pi\sqrt{1-x^2}}dx=\delta_{j,k},\quad k,j\in\mathbb{N}.
\end{equation}

As for uniform boundedness, considering the expressions $C_0(x)=T_0(x)=1$ and $C_n(x)=\sqrt{2}\cos(n\cdot\arccos(x)),n=1,2,\dots$, it is evident that regardless of the highest degree of the polynomials within the normalized first-kind Chebyshev polynomial system, they all exhibit uniform boundedness with a uniform bound $K=\sqrt{2}$.

\subsection{Comparison of Approximation Effects of Sparse Polynomials Utilizing Both Best Approximation Error and $l_1$-minimization Approach Across Various Parameter Choices}
In this subsection, we utilize the first $N+1$ terms $\{C_0(x),\dots,C_{N}(x)\}$ of $\{C_n(x)\}_{n=0}^\infty$ as the orthonormal polynomial basis for our numerical experiments. We initially construct $l_1$-minimization problems to validate the approximation effects of Theorem \ref{thm4} and Theorem \ref{l1Square} respectively.

For Theorem \ref{thm4}, $\Theta=1$ is the most straightforward choice, which leads to the optimization problem (\ref{Infq}) having the constraint $\varepsilon=\sqrt{m}E_N$. By increasing the value of $\Theta$, we can derive different values of $\varepsilon=\Theta\sqrt{m}E_N$. Setting $\varepsilon=0$ corresponds to the interpolation method with the same number of samples. Comparing various $\varepsilon=\Theta\sqrt{m}E_N$ scenarios with the $\varepsilon=0$ scenario allows us to assess whether the function approximation model of Theorem \ref{thm4} outperforms direct interpolation. In the subsequent experiments, adhering to the function approximation process outlined in Theorem \ref{thm4}, we conducted $500$ trials for different original functions $f$ across various maximum polynomial degree $N$ and sampling number $m$. We calculated the relative error $\frac{||f-p^\#||_\infty}{||f||_\infty}$ in the Chebyshev-norm between each approximation polynomial $p^\#$ and $f$, using the average value as a metric for approximation accuracy.

\begin{figure}
    \centering
\includegraphics[height=8.5cm,width=.99\linewidth]{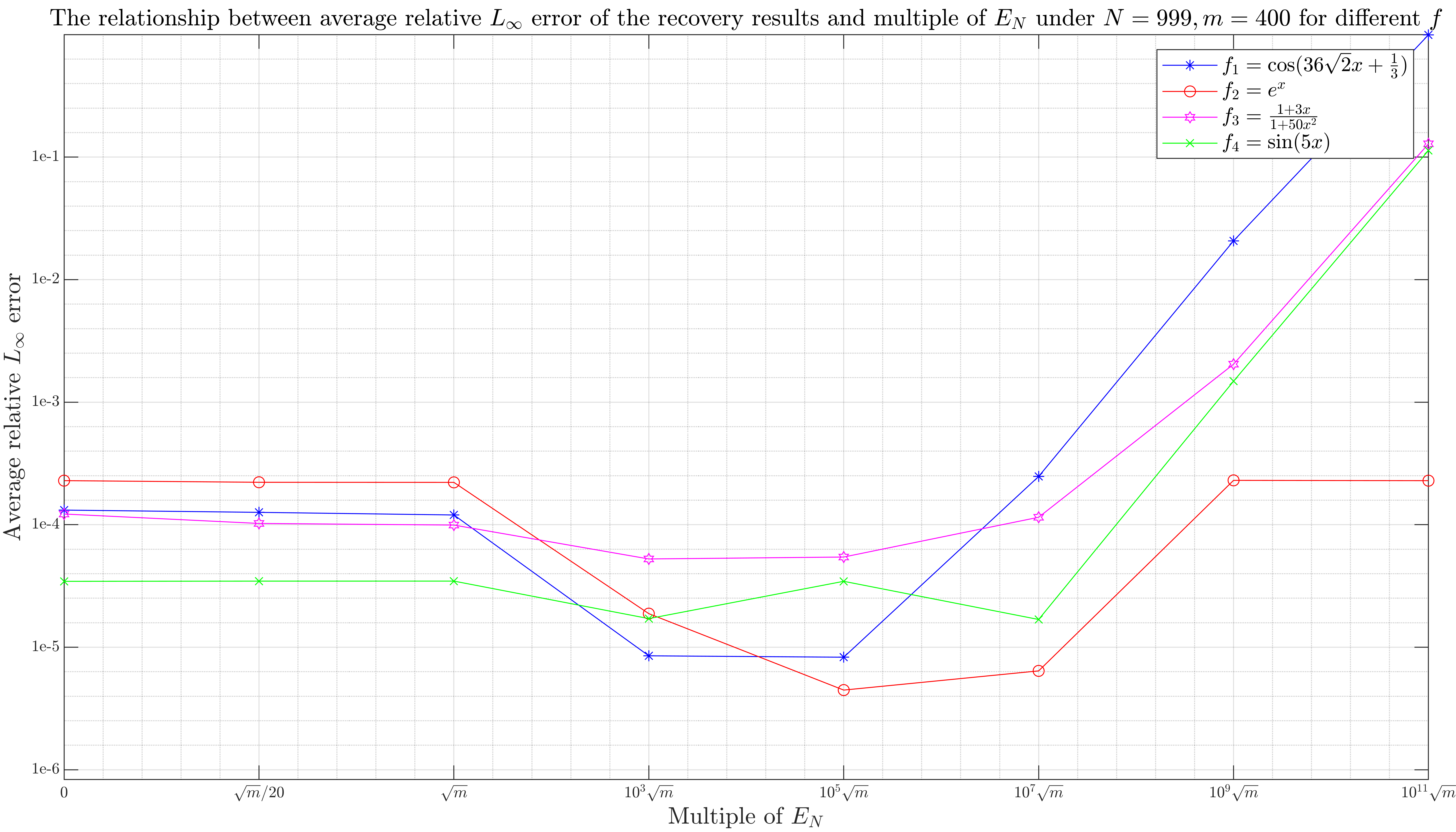}
    \caption{Comparison of the approximation effects for four distinct functions at $N=999,m=400$ utilizing best uniform approximation error and $l_1$-minimization method across various values of $\varepsilon$}
    \label{fig3}
\end{figure}

 We first experiment with increasing the value of $\Theta$ under the same conditions of maximum polynomial degree $N=999$ and sampling number $m=400$, selecting $\Theta=1,10^3,10^5,10^7,10^9,10^{11}$, which correspond to $\varepsilon=\sqrt{m}E_N,10^3\sqrt{m}E_N,10^5\sqrt{m}E_N,10^7\sqrt{m}E_N,10^9\\\sqrt{m}E_N,10^{11}\sqrt{m}E_N$. We performed experiments on different functions and show the results in Figure \ref{fig3}.

Figure \ref{fig3} illustrates how the average Chebyshev-norm relative error changes as the multiple of $E_N$ in the upper bound $\varepsilon$ of the constraint condition for the optimization problem (\ref{Infq}) increases. It is observed that for these specific values of $N$ and $m$, the approximation accuracy to the four different functions initially improves and then declines. Since $N$ is sufficiently large, $E_N$ is typically within the range of $10^{-15}$ to $10^{-14}$. When $\Theta$ is of a moderate size, it does not lead to significant distortion in the results. However, when $\Theta$ becomes excessively large, a noticeable decrease in the approximation accuracy is observed across all four functions. We find that at $\varepsilon=10^5\sqrt{m}E_N$, the average relative errors corresponding to all four functions achieve smaller ones, which is markedly better than the control case where $\varepsilon=0$. This indicates that appropriately increasing the value of $\Theta$ could significantly enhance the approximation accuracy.

\begin{figure}
	\centering
\includegraphics[height=8.5cm,width=.99\textwidth]{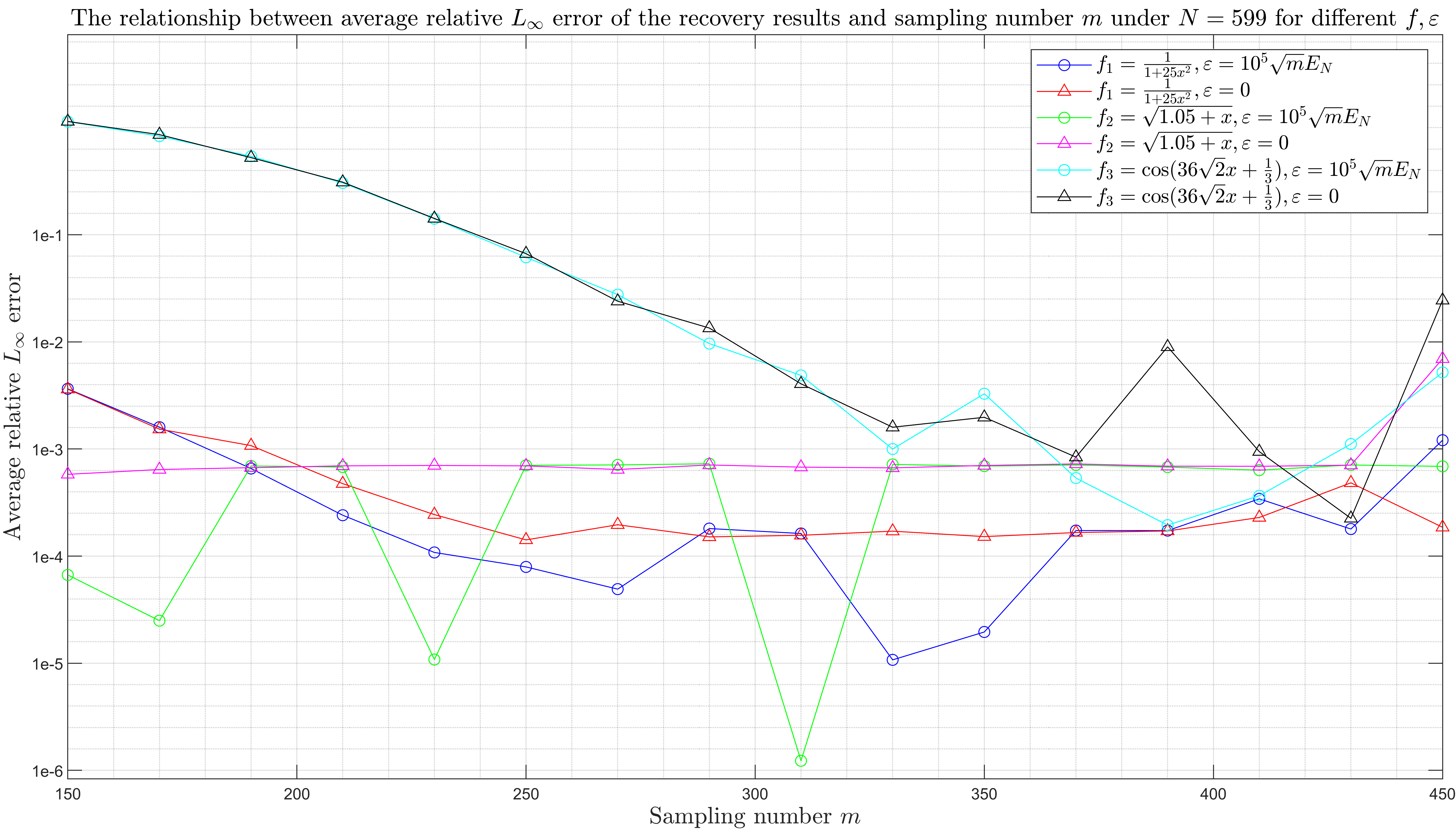}
	\caption{Comparison of the approximation effects for three distinct functions at $N=599$, where $\varepsilon$ is set to $10^5\sqrt{m}E_N$ and $0$ for model (\ref{Infq}), with the sampling number increasing from $150$ to $450$}
	\label{fig4}
\end{figure}

We chose the better-performing $\Theta=10^5$ and, with the maximum polynomial degree $N$ fixed at $599$, observed the comparison of three different functions under the optimization problem (\ref{Infq}) with the upper bound $\varepsilon$ set to $0$ and $10^5\sqrt{m}E_N$ as $m$ increases, as depicted in Figure \ref{fig4}.

While Figure \ref{fig4} indicates that in some instances, the approximation accuracy with $\varepsilon=10^5\sqrt{m}E_N$ does not markedly surpass that of $\varepsilon=0$ for the same combination of $f,N,m$, the former is rarely outperformed. Conversely, varying values of $N$ lead to different $E_N$, making it impractical to set a fixed $\Theta$ that applies universally. As a result, adjustments are necessary based on specific circumstances. Nonetheless, the insights from Figure \ref{fig3} are valuable. By tailoring the upper bound $\varepsilon$ in the optimization problem (\ref{Infq}) to the situation at hand, we can attain approximation results that are considerably more accurate than those obtained through direct interpolation with an equivalent number of samples.

\begin{figure}
	\centering
	\subfigure[$N=599,m=300$]{
		\includegraphics[height=9.5cm,width=.99\linewidth]{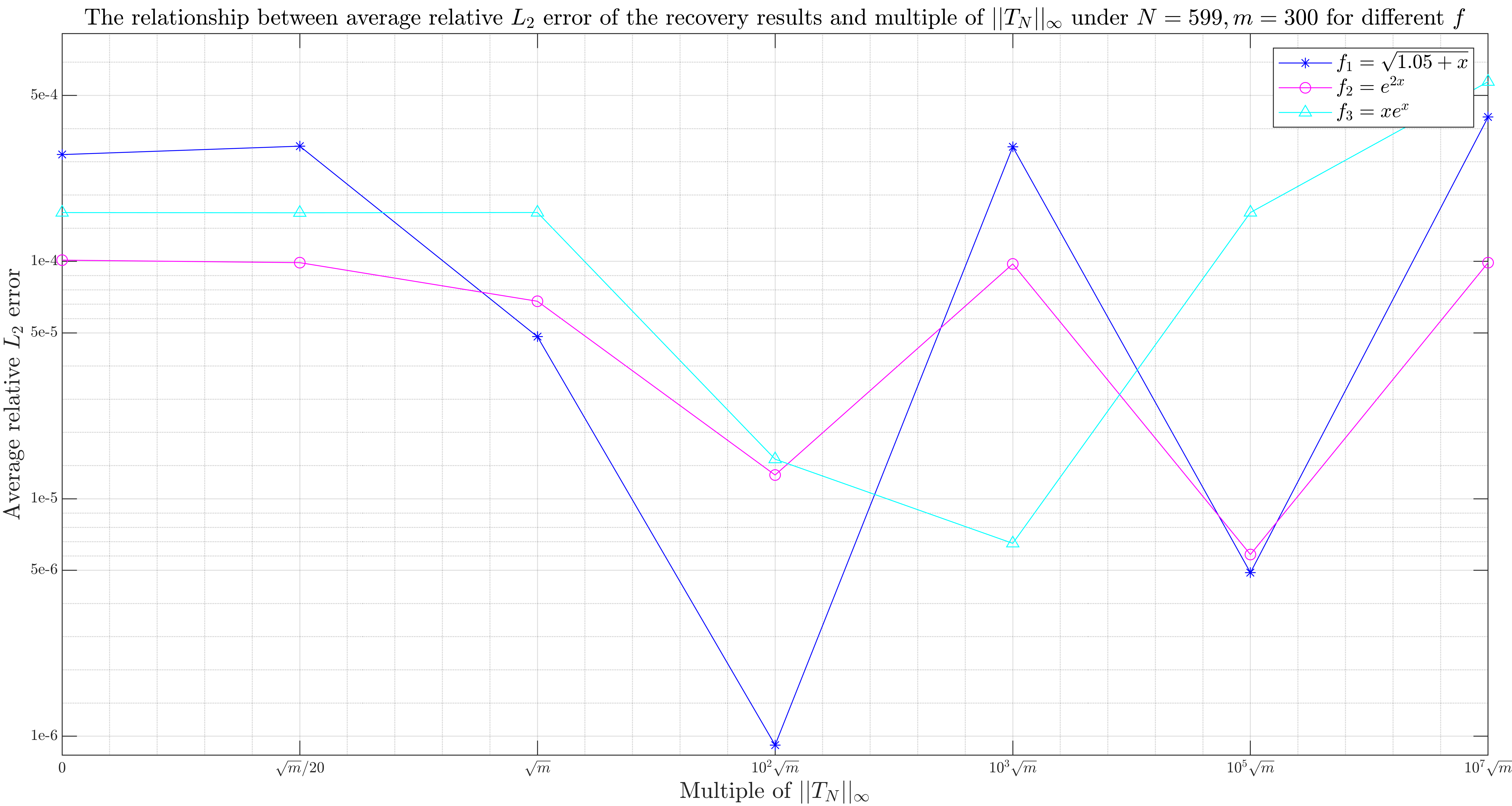}
	}
	\hspace{2em} 
	\subfigure[$N=799,m=400$]{
		\includegraphics[height=9.5cm,width=.99\linewidth]{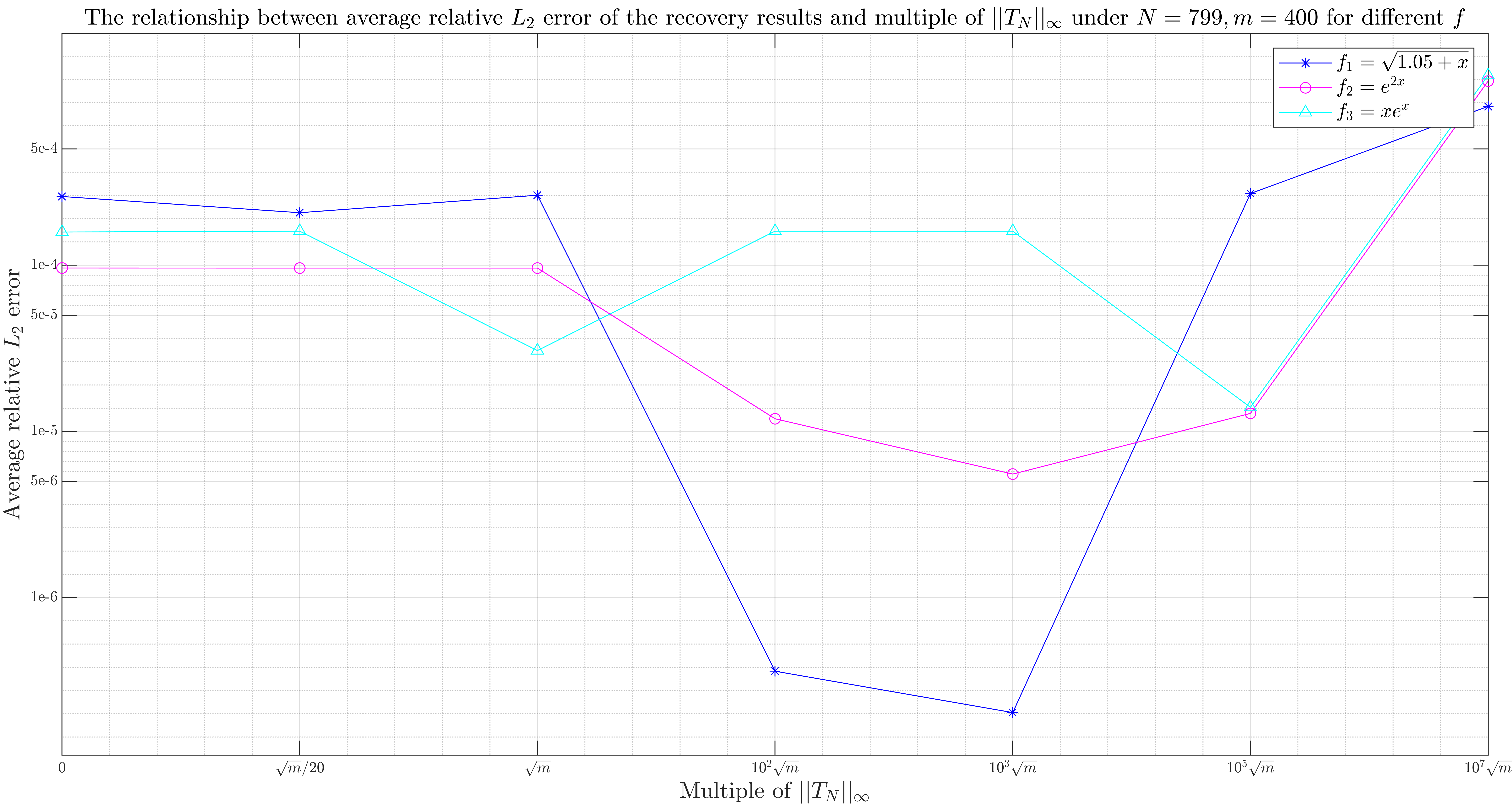}
	}
	\caption{Comparison of the approximation effects for three distinct functions across two different sets of $N,m$ parameters, utilizing the best square approximation error and the $l_1$-minimization method with varying values of $\varepsilon$}
	\label{fig5}
\end{figure}

For the numerical validation of Theorem \ref{l1Square}, we persist in employing $\{C_0(x),\dots,C_{N}(x)\}$ as the orthonormal polynomial basis. For each set of $f,N,m$, we executed $500$ experiments, calculating the relative error $\frac{||f-p^\#||_2}{||f||_2}$ in the $L_2$-norm between each approximation result $p^\#$ and $f$, and considering the average value as an indicator of approximation accuracy. Guided by the numerical experiments based on the $l_1$-minimization method with the best uniform approximation error, we explore various values of $\Theta$ in Theorem \ref{l1Square}, corresponding to various upper bounds $\varepsilon=\Theta\sqrt{m}||T_N||_\infty$ on the right side of the optimization problem (\ref{Squareq}), and compared the approximation accuracy in the $L_2$-norm, as depicted in Figure \ref{fig5}.

\begin{figure}
	\centering
	\subfigure[$N=399$]{
		\includegraphics[height=9.5cm,width=.99\linewidth]{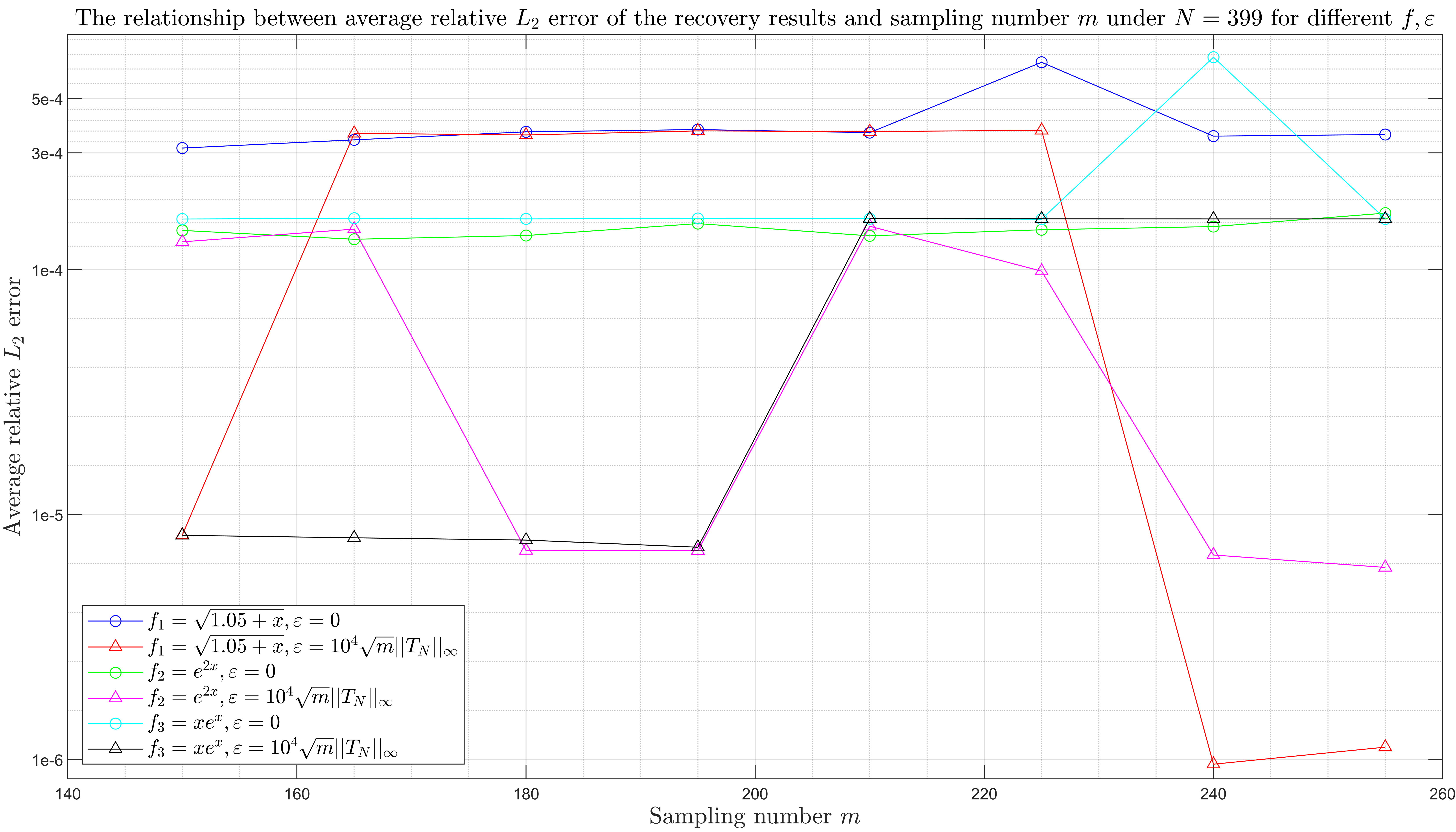}
	}
	\hspace{2em} 
	\subfigure[$N=599$]{
		\includegraphics[height=9.5cm,width=.99\linewidth]{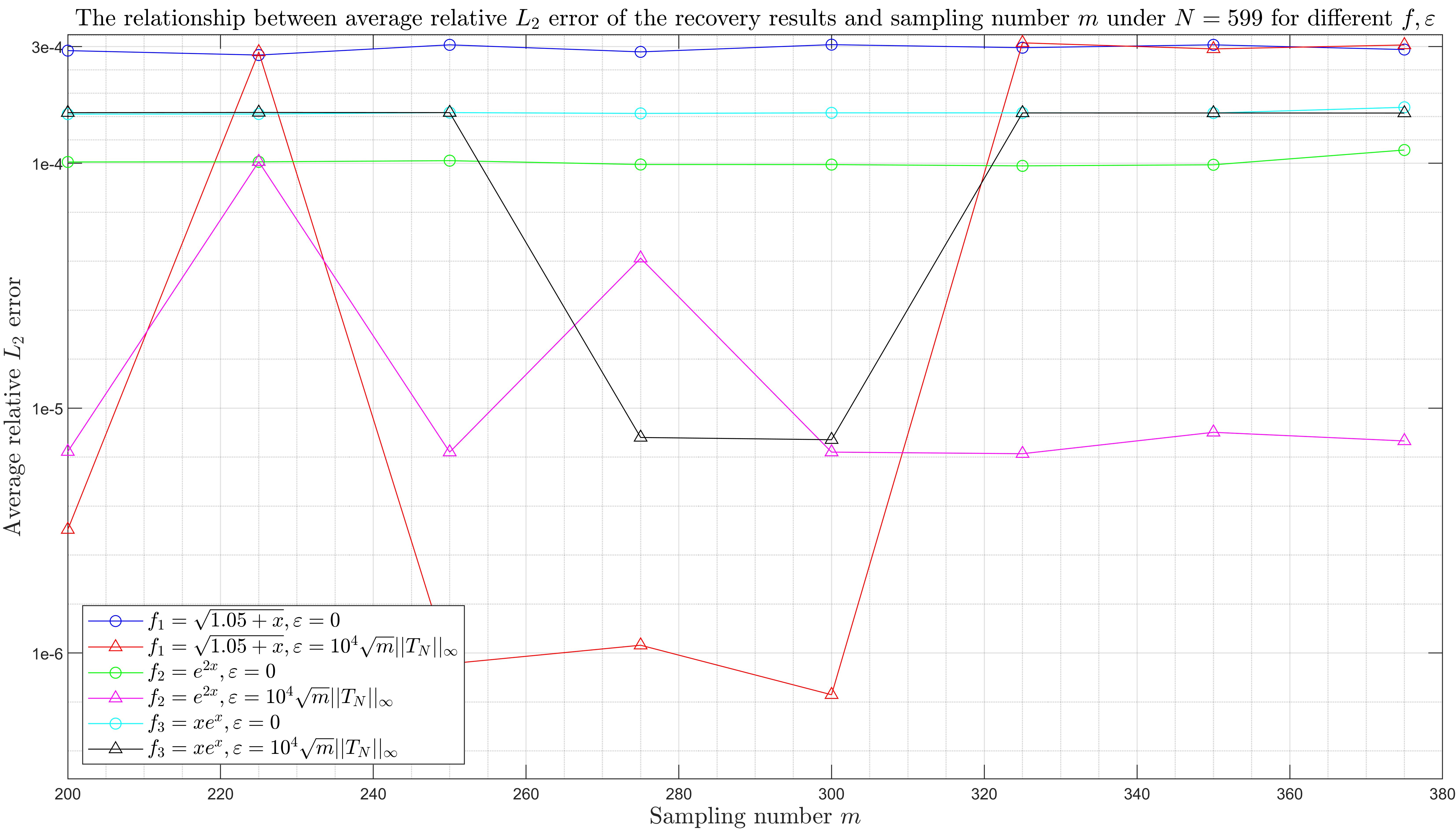}
	}
	\caption{Comparison of the approximation effects for three distinct functions at $N=399$ and $N=599$, utilizing the best square approximation error and the $l_1$-minimization method, with $\varepsilon$ set to $10^4\sqrt{m}||T_N||_\infty$ and $0$, as the sampling number $m$ increases}
	\label{fig6}
\end{figure}

As observed in Figure \ref{fig5}, appropriately increasing $\Theta$ in the approximation method outlined in Theorem \ref{l1Square} can significantly enhance the approximation accuracy. Concentrating on the same three functions from the experiments depicted in Figure \ref{fig5}, we set $N$ to $400$ and $600$, and allowed the sampling number $m$ to vary within an appropriate range. We then compared the average $L_2$-norm relative errors under the conditions $\varepsilon=10^4\sqrt{m}||T_N||_\infty$ and the control scenario where $\varepsilon=0$, showing the results in Figure \ref{fig6}.

The results in Figure \ref{fig6} demonstrate that for different maximum polynomial degrees $N$ and various sampling numbers $m$, the approximation accuracy is generally higher when the upper bound of the optimization problem (\ref{Squareq}) is set to $\varepsilon=10^4\sqrt{m}||T_N||_\infty$ compared to when $\varepsilon=0$, under identical conditions of $f,N,m$. This strongly validates that appropriately increasing the value of $\Theta$ can lead to more accurate approximation results.

Synthesizing the outcomes from these two experimental sets, we conclude that the function approximation methods based on best approximation error and $l_1$-minimization, as presented in Theorem \ref{thm4} and Theorem \ref{l1Square}, are indeed effective and maintain a certain level of precision. Employing the standard $\Theta=1$ in both theorems yields approximation results that are not significantly different from those of interpolation with the same sampling number. However, by moderately increasing the value of $\Theta$, both approximation methods can achieve a marked enhancement in accuracy across a range of $f,N,m$ configurations. Additionally, the choice of $\Theta$ should be cautiously made. It should not be excessively large initially to avoid potential divergence in solving the $l_1$-minimization problem. Moreover, $\Theta$ should be selected based on the specific magnitudes of $E_N$ or $||T_N||_\infty$ , as there is no universal value of $\Theta$ that fits all scenarios optimally.

\subsection{Comparison of Approximation Effects of Sparse Polynomials Utilizing Both Best Approximation Error and Weighted $l_1$-minimization Approach Across Different Weights and Various Parameter Choices}
Moving forward, we focus on conducting numerical experiments to validate the two theorems that are based on solving the weighted $l_1$-minimization problem: Theorem \ref{thm4}'s variant, Theorem \ref{l1Infweight}, and Theorem \ref{l1Square}'s variant, Theorem \ref{l1Squareweight}. We continue to utilize the normalized first-kind Chebyshev polynomial set $\{C_0(x),\dots,C_{N}(x)\}$ on interval $[-1,1]$ as our basis.

Starting with the evaluation of the function approximation performance of Theorem \ref{l1Infweight}, we consider two weight vectors: $\boldsymbol{\omega}^1$ where $\boldsymbol{\omega}^1_i=\sqrt{i}$ for $i=1,2,\dots,N+1$, and $\boldsymbol{\omega}^2$ where $\boldsymbol{\omega}^2_i=\frac{i+1}{2}$ for $i=1,2,\dots,N+1$. With our selected basis, these weight selection methods satisfy the condition $\omega_j\geqslant ||\psi_j||_\infty$ as required by Theorem \ref{l1Infweight} and Theorem \ref{l1Squareweight}, so we also consider these two weights when testing Theorem \ref{l1Squareweight}. For each set consisting of the approximated function $f$, the maximum polynomial degree $N$, and the sampling number $m$, we calculate the average relative error $\frac{||f-p^\#||_\infty}{||f||_\infty}$ of the approximation polynomial $p^\#$. In a single independent repeated experiment, we deem an approximation successful if the relative error is less than $5e^{-4}$. Each set of experiments is repeated $500$ times to determine the average value and to compute the success rate.

\begin{table}[H]
	\centering
	\caption{The results of approximating different functions using the weighted $l_1$-minimization method, based on the best uniform approximation error, under two weights at $N=599$ and $m=300$}
	\label{tab3}
	\begin{tabular}{ccccc}
		\toprule
		$f$ (The approximated function) &$\Theta$& Weight vector & Average relative error &Success rate \\
		\midrule
		$f=\sqrt{1.05+x}$ & 1 & Unweighted &$6.9795\times 10^{-4}$ &0.24\\
		$f=\sqrt{1.05+x}$ & $10^5$& Unweighted &$7.0390\times 10^{-4}$ & 0.224\\
		$f=\sqrt{1.05+x}$ & 1& $\boldsymbol{\omega}^1$& $2.0630\times 10^{-3}$&0\\
		$f=\sqrt{1.05+x}$ & $10^5$&$\boldsymbol{\omega}^1$& \textcolor{red}{$2.0189\times 10^{-6}$}&\textcolor{red}{1}  \\
		$f=\sqrt{1.05+x}$ & 1 &$\boldsymbol{\omega}^2$& $3.0368\times 10^{-3}$&0\\
        $f=\sqrt{1.05+x}$ & $10^5$&$\boldsymbol{\omega}^2$& $3.3784\times 10^{-3}$& 0 \\
		$f=\log(\sin(10x)+2)+\sin(x)$ & 1 & Unweighted & $3.0574\times 10^{-4}$&0.93\\
		$f=\log(\sin(10x)+2)+\sin(x)$ & $10^5$& Unweighted &$8.0949\times 10^{-5}$ & 0.962\\
		$f=\log(\sin(10x)+2)+\sin(x)$ & 1& $\boldsymbol{\omega}^1$&$2.5888\times 10^{-3}$ &0.002\\
		$f=\log(\sin(10x)+2)+\sin(x)$ & $10^5$&$\boldsymbol{\omega}^1$& \textcolor{red}{$1.7109\times 10^{-5}$}& \textcolor{red}{1} \\
		$f=\log(\sin(10x)+2)+\sin(x)$ & 1 &$\boldsymbol{\omega}^2$&$4.4160\times 10^{-3}$ &0\\
        $f=\log(\sin(10x)+2)+\sin(x)$ & $10^5$&$\boldsymbol{\omega}^2$& \textcolor{red}{$2.9866\times 10^{-5}$}& \textcolor{red}{0.996} \\
		\bottomrule
	\end{tabular}
\end{table}

\begin{table}[H]
	\centering
	\caption{The results of approximating different functions using the weighted $l_1$-minimization method, based on the best uniform approximation error, under two weights at $N=799$ and $m=400$}
	\label{tab4}
	\begin{tabular}{ccccc}
		\toprule
		$f$ (The approximated function) &$\Theta$&Weight vector& Average relative error &Success rate \\
		\midrule
		$f=\frac{1}{1+25x^2}$ & 1 & Unweighted & $1.4442\times 10^{-4}$&1\\
		$f=\frac{1}{1+25x^2}$ & $10^5$& Unweighted & $4.9006\times 10^{-6}$& 1\\
		$f=\frac{1}{1+25x^2}$ & 1& $\boldsymbol{\omega}^1$& $9.1198\times 10^{-4}$&0.014\\
		$f=\frac{1}{1+25x^2}$ & $10^5$&$\boldsymbol{\omega}^1$& \textcolor{red}{$7.8976\times 10^{-7}$}& \textcolor{red}{1} \\
		 $f=\frac{1}{1+25x^2}$ & 1 &$\boldsymbol{\omega}^2$&$2.3371\times 10^{-3}$ &0\\
        $f=\frac{1}{1+25x^2}$ & $10^5$&$\boldsymbol{\omega}^2$& $2.5295\times 10^{-3}$& 0 \\
		$f=\cos(36\sqrt{2}x+\frac{1}{3})$ & 1 & Unweighted & $1.9258\times 10^{-4}$&0.974\\
		$f=\cos(36\sqrt{2}x+\frac{1}{3})$ & $10^5$& Unweighted & $3.0894\times 10^{-5}$& 0.998\\
		$f=\cos(36\sqrt{2}x+\frac{1}{3})$ & 1& $\boldsymbol{\omega}^1$&$2.0015\times 10^{-3}$ &0.032\\
		$f=\cos(36\sqrt{2}x+\frac{1}{3})$ & $10^5$&$\boldsymbol{\omega}^1$& \textcolor{red}{$3.8657\times 10^{-6}$}& \textcolor{red}{0.998}\\
		$f=\cos(36\sqrt{2}x+\frac{1}{3})$ & 1 &$\boldsymbol{\omega}^2$& \textcolor{red}{$9.5884\times 10^{-6}$}&\textcolor{red}{1}\\
        $f=\cos(36\sqrt{2}x+\frac{1}{3})$& $10^5$&$\boldsymbol{\omega}^2$&\textcolor{red}{$9.2328\times 10^{-6}$} & \textcolor{red}{1}\\
		\bottomrule
	\end{tabular}
\end{table}

With the same approximated function $f$ and $\Theta$, we have highlighted in red the instances where the average relative error and success rate of the weighted situations outperform the unweighted situations. The results in Table \ref{tab3} demonstrate that after increasing $\Theta$ to $10^5$, which showed better performance in the previous experiments, the introduction of appropriate weights may lead to a notable enhancement in approximation accuracy. Additionally, we have carried out tests under a different set of $N,m$ parameters, as illustrated in Table \ref{tab4}.

The results in Table \ref{tab4} indicate that with $\Theta=10^5$, the introduction of weights continues to produce excellent outcomes. By synthesizing the findings from Table \ref{tab3} and Table \ref{tab4}, we deduce that properly elevating the value of $\Theta$ in Theorem \ref{l1Infweight} and employing weights that align with the stipulations of Theorem \ref{l1Infweight} for weighted $l_1$-minimization frequently results in enhanced approximation accuracy.

\begin{table}[H]
	\centering
	\caption{The results of approximating different functions using the weighted $l_1$-minimization method, based on the best square approximation error, under two weights at $N=799$ and $m=400$}
	\label{tab5}
	\begin{tabular}{ccccc}
		\toprule
		$f$ (The approximated function) &$\Theta$&Weight vector& Average relative error &Success rate \\
		\midrule
		$f=\frac{1}{1+25x^2}$ & 1 & Unweighted &$9.9227\times 10^{-5}$ &1\\
		$f=\frac{1}{1+25x^2}$ & $10^5$& Unweighted & $4.8450\times 10^{-6}$&1 \\
		$f=\frac{1}{1+25x^2}$ & 1& $\boldsymbol{\omega}^1$&$2.2710\times 10^{-4}$ & 0.972\\
		$f=\frac{1}{1+25x^2}$ & $10^5$&$\boldsymbol{\omega}^1$&\textcolor{red}{$1.7887\times 10^{-6}$} &\textcolor{red}{1} \\
		$f=\frac{1}{1+25x^2}$ & 1 &$\boldsymbol{\omega}^2$& $5.5637\times 10^{-4}$&0.434\\
        $f=\frac{1}{1+25x^2}$ & $10^5$&$\boldsymbol{\omega}^2$&\textcolor{red}{$3.9910\times 10^{-6}$} &0.998  \\
		$f=\cos(36\sqrt{2}x+\frac{1}{3})$ & 1 & Unweighted & $5.1107\times 10^{-5}$& 0.998\\
		$f=\cos(36\sqrt{2}x+\frac{1}{3})$ & $10^5$& Unweighted &$8.0470\times 10^{-5}$ & 0.998\\
		$f=\cos(36\sqrt{2}x+\frac{1}{3})$ & 1& $\boldsymbol{\omega}^1$& $8.8177\times 10^{-4}$& 0.382\\
		$f=\cos(36\sqrt{2}x+\frac{1}{3})$ & $10^5$&$\boldsymbol{\omega}^1$&\textcolor{red}{$1.4300\times 10^{-5}$} & \textcolor{red}{1}\\
		$f=\cos(36\sqrt{2}x+\frac{1}{3})$ & 1 &$\boldsymbol{\omega}^2$&\textcolor{red}{$2.8717\times 10^{-6}$} &\textcolor{red}{1}\\
        $f=\cos(36\sqrt{2}x+\frac{1}{3})$& $10^5$&$\boldsymbol{\omega}^2$&\textcolor{red}{$3.1737\times 10^{-6}$} & \textcolor{red}{1}\\
		$f=\sqrt{1.05+x}$ & 1 & Unweighted &$2.6360\times 10^{-4}$ &0.998\\
		$f=\sqrt{1.05+x}$ & $10^5$& Unweighted &$2.6991\times 10^{-4}$ &1 \\
		$f=\sqrt{1.05+x}$ & 1& $\boldsymbol{\omega}^1$& \textcolor{red}{$3.5372\times 10^{-5}$}&\textcolor{red}{1} \\
		$f=\sqrt{1.05+x}$ & $10^5$&$\boldsymbol{\omega}^1$&$8.2742\times 10^{-4}$ & 0.022\\
		$f=\sqrt{1.05+x}$ & 1 &$\boldsymbol{\omega}^2$&\textcolor{red}{$5.2161\times 10^{-5}$} &\textcolor{red}{1}\\
        $f=\sqrt{1.05+x}$ & $10^5$&$\boldsymbol{\omega}^2$&$1.3122\times 10^{-3}$ &0  \\
		$f=\log(\sin(10x)+2)+\sin(x)$ & 1 & Unweighted & $2.6140\times 10^{-5}$& 1\\
		$f=\log(\sin(10x)+2)+\sin(x)$ & $10^5$& Unweighted & $1.4076\times 10^{-5}$& 1\\
		$f=\log(\sin(10x)+2)+\sin(x)$ & 1& $\boldsymbol{\omega}^1$&$1.3879\times 10^{-3}$ & 0.088\\
		$f=\log(\sin(10x)+2)+\sin(x)$ & $10^5$&$\boldsymbol{\omega}^1$&\textcolor{red}{$1.2234\times 10^{-5}$} & \textcolor{red}{1}\\
		$f=\log(\sin(10x)+2)+\sin(x)$ & 1 &$\boldsymbol{\omega}^2$&$2.6167\times 10^{-3}$ & 0\\
        $f=\log(\sin(10x)+2)+\sin(x)$& $10^5$&$\boldsymbol{\omega}^2$&\textcolor{red}{$1.3211\times 10^{-5}$} & \textcolor{red}{1}\\
		
		\bottomrule
	\end{tabular}
\end{table}

Proceeding to assess the function approximation performance of Theorem \ref{l1Squareweight}, we persist in utilizing weights $\boldsymbol{\omega}^1$ and $\boldsymbol{\omega}^2$. The relative error is expressed as $\frac{||f-p^\#||_2}{||f||_2}$, and we deem an approximation successful if it falls below $5e^{-4}$. Each experimental setup is repeated $500$ times to determine the average relative error and to tally the success rate. Concentrating on the parameters $N=799$ and $m=400$ that showed superior results in the weighted experiments associated with the best uniform approximation, we show the outcomes in Table \ref{tab5}.

In Table \ref{tab5}, $\Theta$ denotes the predetermined positive integer as specified in Theorem \ref{l1Squareweight}. It is observed that setting $\Theta=10^5$ often leads to a reduction in the average relative error and an increase in the success rate when weights are introduced. However, there are instances where weights can have an adverse effect. The performance of the two weight vectors,  $\boldsymbol{\omega}^1$ and $\boldsymbol{\omega}^2$, varies, which supports the notion that different weights exhibit varying levels of effectiveness under different conditions.

When comparing weighted scenarios with unweighted scenarios, we find that introducing weights when appropriately increasing the predetermined parameter $\Theta$ tends to yield more precise approximation outcomes. The experimental results from the two novel weighted methods demonstrate that the performance of different weights is influenced by various factors, including the original function $f$, the maximum polynomial degree $N$, and the sampling number $m$. Consequently, the selection of weights should be tailored to the specific context and no one-size-fits-all weight selection strategy has been identified that is universally applicable across all scenarios.

\section{Summary}
This paper employs best approximation and $l_1$-minimization or weighted $l_1$-minimization methods for sparse polynomial approximation of continuous functions, offering approximation error estimates. While the conclusions presented herein stem from the context of univariate function approximation, their applicability extends to multivariate scenarios. However, estimating the magnitude of error within a multivariate setting encounters challenges, as there currently exists no viable operational algorithm to assess the best uniform approximation error $E_N$ for multivariate continuous functions relative to finite multivariate polynomial spaces.

Additionally, the question of how to develop more precise weight selection strategies that confer a distinct advantage to weighted $l_1$-minimization methods over their unweighted counterparts remains to be resolved. It is possible that theoretical derivations could yield a class of weight selection strategies applicable across various contexts.

\section*{Acknowledgments}
The work is supported by National Natural Science Foundation of China no.12071019.


\end{CJK}
\end{document}